\documentclass[journal]{IEEEtran}
\newtheorem{theorem}{Theorem}
\newtheorem{remark}{Remark}

\usepackage{cite}
\ifCLASSINFOpdf
\else
\fi
\usepackage[cmex10]{amsmath}
\usepackage{algorithmicx}
\usepackage{algpseudocode}
\usepackage{threeparttable}

\usepackage{graphicx}
\begin{document}
%
% paper title
% can use linebreaks \\ within to get better formatting as desired
\title{Optimal Design of Water Distribution Networks by Discrete State Transition Algorithm}
%
%
% author names and IEEE memberships
% note positions of commas and nonbreaking spaces ( ~ ) LaTeX will not break
% a structure at a ~ so this keeps an author's name from being broken across
% two lines.
% use \thanks{} to gain access to the first footnote area
% a separate \thanks must be used for each paragraph as LaTeX2e's \thanks
% was not built to handle multiple paragraphs
%

\author{Xiaojun Zhou
        David Yang Gao,
        and Angus R. Simpson% <-this % stops a space

\thanks{Xiaojun Zhou is with the School of Science, Information Technology and Engineering, University of Ballarat, Victoria 3353, Australia
and School of Information Science and Engineering, Central South University, Changsha 410083, China (tiezhongyu2010@gmail.com)}% <-this % stops a space
\thanks{David Yang Gao is with the School of Science, Information Technology and Engineering, University of Ballarat, Victoria 3353, Australia (d.gao@ballarat.edu.au)}% <-this % stops a space
\thanks{Angus R. Simpson is with the Centre for Applied Modeling in Water Engineering, School of Civil and Environmental Engineering, University of Adelaide, Adelaide,
SA 5005, Australia (angus.simpson@adelaide.edu.au)}% <-this % stops a space
\thanks{Manuscript received ** , 2013; revised **, 2013.}}

% note the % following the last \IEEEmembership and also \thanks -
% these prevent an unwanted space from occurring between the last author name
% and the end of the author line. i.e., if you had this:
%
% \author{....lastname \thanks{...} \thanks{...} }
%                     ^------------^------------^----Do not want these spaces!
%
% a space would be appended to the last name and could cause every name on that
% line to be shifted left slightly. This is one of those "LaTeX things". For
% instance, "\textbf{A} \textbf{B}" will typeset as "A B" not "AB". To get
% "AB" then you have to do: "\textbf{A}\textbf{B}"
% \thanks is no different in this regard, so shield the last } of each \thanks
% that ends a line with a % and do not let a space in before the next \thanks.
% Spaces after \IEEEmembership other than the last one are OK (and needed) as
% you are supposed to have spaces between the names. For what it is worth,
% this is a minor point as most people would not even notice if the said evil
% space somehow managed to creep in.

% The paper headers
% IEEE Transactions on Systems, Man and Cybernetics
\markboth{Journal of \LaTeX\ Class Files, ~Vol.~, No.~, January~2013}%
{Shell \MakeLowercase{\textit{et al.}}: Bare Demo of IEEEtran.cls for Journals}
% The only time the second header will appear is for the odd numbered pages
% after the title page when using the twoside option.
%
% *** Note that you probably will NOT want to include the author's ***
% *** name in the headers of peer review papers.                   ***
% You can use \ifCLASSOPTIONpeerreview for conditional compilation here if
% you desire.

% If you want to put a publisher's ID mark on the page you can do it like
% this:
%\IEEEpubid{0000--0000/00\$00.00~\copyright~2007 IEEE}
% Remember, if you use this you must call \IEEEpubidadjcol in the second
% column for its text to clear the IEEEpubid mark.

% use for special paper notices
%\IEEEspecialpapernotice{(Invited Paper)}

% make the title area
\maketitle

\begin{abstract}
%\boldmath
Optimal design of water distribution networks, which are governed by a series of  linear and nonlinear equations, has been extensively studied in the past decades. Due to their NP-hardness, methods to solve the optimization problem have changed from traditional
mathematical programming to modern intelligent optimization techniques.
In this study, with respect to the model formulation, we have demonstrated that the network system can be reduced to the dimensionality of the number of closed simple loops or required independent paths, and the reduced nonlinear system can be solved efficiently by the Newton-Raphson method.
Regarding the optimization technique, a discrete state transition algorithm (STA) is introduced to solve several cases of water distribution networks.
In discrete STA, there exist four basic intelligent operators, namely, swap, shift, symmetry and substitute as well as the ``risk and restore in probability" strategy.
Firstly, we focus on a parametric study of the restore probability $p_1$ and risk probability $p_2$. To effectively deal with the head pressure constraints, we then investigate the
effect of penalty coefficient and search enforcement on the performance of the algorithm. Based on the experience gained from the training of the Two-Loop network problem, the discrete STA
has successfully achieved the best known solutions for the Hanoi and New York problems. A detailed comparison of our results with those gained by other algorithms is also presented.
\end{abstract}
% IEEEtran.cls defaults to using nonbold math in the Abstract.
% This preserves the distinction between vectors and scalars. However,
% if the journal you are submitting to favors bold math in the abstract,
% then you can use LaTeX's standard command \boldmath at the very start
% of the abstract to achieve this. Many IEEE journals frown on math
% in the abstract anyway.

% Note that keywords are not normally used for peerreview papers.

% For peer review papers, you can put extra information on the cover
% page as needed:
% \ifCLASSOPTIONpeerreview
% \begin{center} \bfseries EDICS Category: 3-BBND \end{center}
% \fi
%
% For peerreview papers, this IEEEtran command inserts a page break and
% creates the second title. It will be ignored for other modes.
\IEEEpeerreviewmaketitle

\section{Introduction}
% The very first letter is a 2 line initial drop letter followed
% by the rest of the first word in caps.
%
% form to use if the first word consists of a single letter:
% \IEEEPARstart{A}{demo} file is ....
%
% form to use if you need the single drop letter followed by
% normal text (unknown if ever used by IEEE):
% \IEEEPARstart{A}{}demo file is ....
%
% Some journals put the first two words in caps:
% \IEEEPARstart{T}{his demo} file is ....
%
% Here we have the typical use of a "T" for an initial drop letter
% and "HIS" in caps to complete the first word.
\IEEEPARstart{P}{ipes}, hydraulic devices (pumps, valves, etc.) and reservoirs are connected in a water distribution network in a complex manner. The physical behavior of a looped network is governed by a set of linear and nonlinear equations, including continuity and energy equations, and head loss functions. The overall planning tasks to be performed in water distribution networks consists of three kinds of problems: layout, design and operation. Although these problems are not independent with each other, they can be formulated and solved separately from a technical point of view since each one can be considered as a parameter when others are being solved.  In this work, we focus on the optimal design problem.\\
\indent Optimal selection of pipe diameters to constitute a water distribution network respecting certain pressure requirements has been shown to be an NP-hard problem \cite{yates1984}, mainly due to two reasons: nonlinear equations and discrete-valued diameters.
A terribly clumsy method for designing pipe network is by enumeration or complete trial and error \cite{gessler1985}.
Traditional methods are to linearize and relax the problem firstly to facilitate the use of linear programming and nonlinear programming, and then they have to round off the solution to the nearest discrete diameters \cite{alperovits1977, morgan1985, goulter1986, kessler1989, eiger1994}. Such algorithms can not guarantee global optima and sometimes cause infeasible solutions. In the last few decades, intelligent optimization techniques including: genetic algorithm \cite{dandy1996, savic1997, simpson1994}, simulated annealing \cite{cunha1999}, shuffled complex evolution \cite{liong2004}, ant colony optimization \cite{zecchin2005,zecchin2006}, harmony search \cite{geem2006}, particle swarm optimization \cite{montalvo2008}, differential evolution \cite{vasan2010} and some of their hybrids \cite{cisty2010,haghighi2011}, have found wide applications in this field.
The advantages of using these stochastic algorithms are: (1) simple representation of a discrete-valued solution; (2) independent to the problem structure to some extent; (3) easy computation due to the only use of the information about the objective function; (4) high probability to gain the global optimum or approximate global optimum in a reasonable amount of time.\\
\indent We introduce the recently developed intelligent optimization algorithm, state transition algorithm (STA) \cite{xzhou2011a,xzhou2011b,xzhou2012}, which shows fantastic performance in continuous function optimization.
In \cite{yang2012}, a discrete STA was proposed to solve the traveling salesman problem, and the results demonstrated that it consumed much less time and had better search ability than
the well-known simulated annealing and ant colony optimization. The goal of this paper is to apply the discrete STA to the optimal design problem of the water distribution networks.\\
\indent This paper is organized as follows. In Section II, the optimization model of water distribution networks is established, including the objective function, decision variables and some constraints.
In Section III, the basic key elements in discrete STA are introduced. It focuses on the intelligent operators of discrete STA and a parametric study of the
``restore probability" and ``risk probability" is the emphasis. How to deal with the constraints and the implementation of the discrete STA for the optimal design problem are illustrated
in Section IV. In Section V, several case studies are given. The Two-Loop network is mainly studied to investigate the effect of penalty coefficient and search enforcement on the performance of the discrete STA. The gained experience is applied to other cases and the results achieved by the proposed discrete STA with other optimization algorithms are presented as well.  Conclusion is derived in Section VI.

% You must have at least 2 lines in the paragraph with the drop letter
% (should never be an issue)
\section{Optimization model formulation of water distribution networks}
\indent For a given layout of pipes and a set of specified demand patterns at the nodes, the optimal design of a water distribution network is to find the combination of commercial pipe sizes which gives the minimum cost, subject to the following constraints:
\begin{itemize}
  \item continuity of flow;
  \item head loss;
  \item conservation of energy;
  \item minimum pressure head.
\end{itemize}
\subsection{The objective function}
\indent Considering that the pipe layout, connectivity and imposed minimum head constraints are known, in the optimal design problem of the water distribution network, the pipe diameters are the only decision variables. As a result, the objective function is assured to be a cost function of pipe diameters
\begin{eqnarray}
 \min_{D_j \in \Omega}~~ f_{obj} = \sum_{j=1}^{NP} L_j c(D_j),
\end{eqnarray}
where, $\Omega$ is a set of commercial pipe sizes, $NP$ is the number of pipes, and $L_j$ is the length of pipe \textit{j}, which is known in this study. $c(D_j)$ indicates that for every commercial pipe size, there is a corresponding cost per unit associated with it.
\subsection{Continuity equation}
Conservation of mass at nodes or junctions in a water distribution network yields a set of linear algebraic equations in terms of flows. At each node, flow continuity should be satisfied,
\begin{eqnarray}
 -\sum Q_{in} + \sum Q_{out} + DM = 0,
\end{eqnarray}
where, $DM$ is the demand at the node, $Q_{in}$ and $Q_{out}$ are the flow entering and leaving the node, respectively.
\subsection{Head loss equation}
The head loss in a pipe in the water distribution network can be computed from a number of empirically obtained equations. The two commonly used equations are the Darcy-Weisbach head loss equation and the Hazen-Williams head loss equation. The general expression for the head loss in a pipe \textit{j} located between nodes \textit{i} and \textit{k} is given by
\begin{eqnarray}
H_i - H_k = r_j Q_j|Q_j|^{\alpha-1} = \omega \frac{L_j}{C^{\alpha} D_j^{\beta}}  Q_j|Q_j|^{\alpha-1},
\end{eqnarray}
where, $H_i$ and $H_k$ are nodal pressure head at the end of the pipe at node \textit{i} and \textit{k} respectively; $r_j$ is called resistance factor for the pipe \textit{j}; $Q_j$ is the flow in pipe \textit{j}; $\omega$ is a numerical conversion constant depending on the units used; $L_j$ is the length of pipe \textit{j}; $C$ is the roughness coefficient; $\alpha$ and $\beta$ are coefficients.\\
\indent For International System of Units (SI), $\omega = 10.6744 $ or $\omega = 10.5088 $, $\alpha = 1/0.54=1.852$ and $\beta = 2.63/0.54=4.871$ are employed in this study using the Hazen-Williams formula.
\subsection{Energy equation}
Energy conservation equations around closed simple loops or between fixed head nodes along required independent paths in a network are nonlinear. Upon traversing a closed simple loop or a required independent path, the sum of pipe head losses around the loop or the path must be zero, which can be expressed as
\begin{eqnarray}
\sum_{j \in L_s} \omega \frac{L_j}{C^{\alpha} D_j^{\beta}}  Q_j|Q_j|^{\alpha-1} - \sum_{j \in L_s} EL_j = 0,
\end{eqnarray}
where, $L_s$ is the indices of pipes in a closed simple loop or a required independent path; $EL_j$ is the hydraulic grade line at the reservoir \textit{j}.
\subsection{Minimum pressure head}
The minimum pressure head constraints at each node are given as follows
\begin{eqnarray}
H_i \geq H_{i\min}, \forall i = 1, \cdots, NJ,
\end{eqnarray}
where, $H_{i\min}$ is known, and $NJ$ is the number of nodes.
\section{A brief review of the discrete state transition algorithm}
Let's consider the following unconstrained integer optimization problem
\begin{eqnarray}\label{eqn1}
\min f(x),
\end{eqnarray}
where, $x = (x_1, \cdots, x_n)$, $x_i \in \mathcal{I} \subset \mathcal{Z}^m$, $i = 1, \cdots, n$, and $f(x)$ is a real-valued function.
\subsection{The framework of the discrete state transition algorithm}
\indent If a solution to a specific optimization problem is described as a state, then the transformation to update the solution becomes a state transition. Without loss of generality, the unified form of discrete state transition algorithm can be described as
\begin{eqnarray}
\left \{ \begin{array}{ll}
x_{k+1}= A_{k}(x_{k}) \bigoplus B_{k}(u_{k})\\
y_{k+1}= f(x_{k+1})
\end{array} \right.,
\end{eqnarray}
where, $x_{k} \in \mathcal{Z}^{n}$ stands for a current state, corresponding to a solution of a discrete optimization problem; $u_{k}$ is a function of $x_{k}$ and historical states; $A_{k}(\cdot)$, $B_{k}(\cdot)$ are transformation operators, which are usually state transition matrixes; $\bigoplus$ is a operation, which is admissible to operate on two states; $f$ is the cost function or evaluation function.\\
\indent As a intelligent optimization algorithm, the discrete state transition algorithm have the following five key elements:\\
\indent (1) Representation of a solution. In discrete STA, we choose special representations, that is, the permutation of the set $\{1,2,\cdots,n\}$, which can be easily manipulated by some intelligent operators. The reason that we call the operators ``intelligent" is due to their geometrical property (swap, shift, symmetry and substitute), and a intelligent operator has the same geometrical function for different representations. A big advantage of such representations and operators is that, after each state transformation, the newly created state is always feasible, avoiding the trouble into rounding off a continuous solution in other cases.\\
\indent (2) Sampling in a candidate set. When a transformation operator is exerted on a current state, the next state is not deterministic, that is to say, there are possibly different choices for the next state. It is not difficult to imagine that all possible choices will constitute a candidate set, or a ``neighborhood". Then we execute several times of transformation (called search enforcement $SE$) on current state, to sampling in the ``neighborhood". Sampling is a very important factor in state transition algorithm, which can reduce the search space and avoid enumeration.\\
\indent (3) Local exploitation and global exploration. In optimization algorithms, it is quite significant to design good local and global operators. The local exploitation can guarantee high precision of a solution and convergent performance of a algorithm, and the global exploration can avoid getting trapped into local minima or prevent premature convergence.
In discrete optimization, it is extremely difficult to define a ``good" local optimal solution due to its dependence on a problem's structure, which leads to the same difficulty in the definition of local exploitation and global exploration. Anyway, in the discrete state transition algorithm, we define the slow change to current solution by a transformation as local exploitation, while the big change to current solution by a transformation as global exploration.\\
\indent (4) Self learning and regular communication. State transition algorithm behaves in two styles, one is individual-based, the other is population-based, which is certainly a extended version. The individual-based state transition algorithm focuses on self learning, in other words, with emphasis on the operators' designing and dynamic adjustment (details given in the following). Undoubtedly, communication among different states is a promising strategy for state transition algorithm, as indicated in \cite{xzhou2012}. Through communication, states can share information and cooperate with each other. However, how to communicate and when to communicate are key issues. In continuous state transition algorithm, intermittent exchange strategy was proposed, which means that states communicate with each other at a certain frequency in a regular way.\\
\indent (5) Dynamic adjustment. It is a potentially useful strategy for state transition algorithm. In the iteration process of an intelligent algorithm, the fitness value can decrease sharply in the early stage, but it stagnates in the late stage, due to the static environment. As a result, some perturbation should be added to activate the environment. In fact, dynamic adjustment can be understood and implemented in various ways. For example, the alternative use of different local and global operators is dynamic adjustment to some extent. Then, we can change the search enforcement, vary the cost function, reduce the dimension, etc. Of course, ``risk a bad solution in probability" is another dynamic adjustment, which is widely used in simulated annealing (SA). In SA, the Metropolis criterion \cite{metroplis1953} is used to accept a bad solution:
\begin{eqnarray}
\mathrm{probability}~~ p = \mathrm{exp}(\frac{-\triangle E}{k_B T}),
\end{eqnarray}
where, $\triangle E = f(x_{k+1}) - f(x_{k})$, $k_B$ is the Boltzmann probability  factor, $T$ is the temperature to regulate the process of annealing. In the early stage, temperature is high, and it has big probability to accept a bad solution, while in the late stage, temperature is low, and it has very small probability to accept a bad solution, which is the key point to guarantee the convergence. We can see that the Metropolis criterion has the ability to escape from local optimality, but on the other hand, it will miss some ``good solutions" as well.
\\
\indent In this study, we focus on the individual-based STA, and the main process of discrete STA is shown in the pseudocode as follows
\begin{algorithmic}[1]
\Repeat
    \State {[Best,fBest] $\gets$ swap(fcn,Best,fBest,SE,n,$m_a$)}
    \State {[Best,fBest] $\gets$ shift(fcn,Best,fBest,SE,n,$m_b$)}
    \State {[Best,fBest] $\gets$ symmetry(fcn,Best,fBest,SE,n,$m_c$)}
    \State {[Best,fBest] $\gets$ substitute(fcn,Best,fBest,SE,set,n,$m_d$)}
    \If{ fBest $<$ fBest$^{*}$} \Comment{greedy criterion}
    \State {Best$^{*}$ $\gets$ Best}
    \State {fBest$^{*}$ $\gets$ fBest}
    \EndIf
    \If{$rand < p_1$} \Comment{restore in probability}
    \State {Best $\gets$ Best$^{*}$}
    \State {fBest $\gets$ fBest$^{*}$}
    \EndIf
\Until{the maximum number of iterations is met}
\end{algorithmic}
\quad \\
\indent As for detailed explanations, swap function in above pseudocode is given as follows for example
\begin{algorithmic}[1]
\State{State $\gets$ op\_swap(Best,SE,n,$m_a$)}
\State{[newBest,fnewBest] $\gets$ fitness(funfcn,State)}
\If{fnewBest $<$ fBest}\Comment{greedy criterion}
    \State{Best $\gets$ newBest}
    \State{fBest $\gets$ fnewBest}
\Else
    \If{$rand < p_2$}\Comment{risk in probability}
    \State{Best $\gets$ newBest}
    \State{fBest $\gets$ fnewBest}
    \EndIf
\EndIf
\end{algorithmic}
\indent From the pseudocodes, we can find that in discrete STA, in the whole, ``greedy criterion" is adopted to keep the incumbent ``$\mathrm{Best}^{*}$", in the partial, a bad solution ``$\mathrm{Best}$" is accepted in each inner state transformation at a  probability $p_2$, and in the same while, the ``$\mathrm{Best}^{*}$" is restored in the outer iterative process at another probability $p_1$. The ``risk a bad solution in probability" strategy aims to escape from local optimal, while the ``greedy criterion" and ``restore the incumbent best solution in probability" are to guarantee a good convergence.
\subsection{The representation, local and global operators}
\begin{figure}[!htbp]
  \centering
  % Requires \usepackage{graphicx}
  \includegraphics[width=9cm]{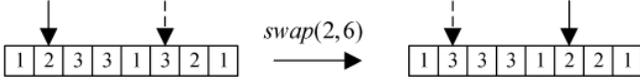}
  \caption{illustration of the swap transformation}
  \label{dvsswap}
\end{figure}

\begin{figure}[!htbp]
  \centering
  % Requires \usepackage{graphicx}
  \includegraphics[width=9cm]{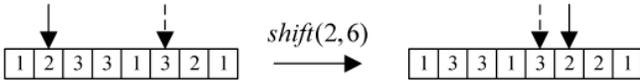}
  \caption{illustration of the shift transformation}\label{dvsshift}
\end{figure}

\begin{figure}[!htbp]
  \centering
  % Requires \usepackage{graphicx}
  \includegraphics[width=9cm]{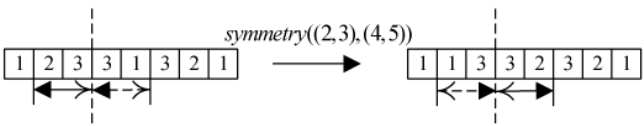}
  \caption{illustration of the symmetry transformation}\label{dvssymmetry}
\end{figure}

\begin{figure}[!htbp]
  \centering
  % Requires \usepackage{graphicx}
  \includegraphics[width=9cm]{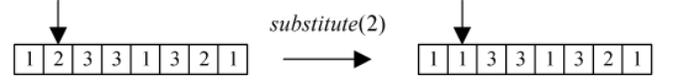}
  \caption{illustration of the substitute transformation}
  \label{dvssubstitute}
\end{figure}

In discrete STA, we use the index of the a commercial size  as a representation for a solution to the optimal design problem. For example, if there are 8 pipes and for each pipe there are 3 choices,
then the details of four special geometric operators are defined as follows
\\
(1) Swap transformation
\begin{eqnarray}
x_{k+1}= A^{swap}_{k}(m_a) x_{k},
\end{eqnarray}
where, $A^{swap}_{k} \in \mathcal{Z}^{n \times n}$ is called swap permutation matrix, $m_a$ is a constant integer called swap factor to control the maximum number of positions to be exchanged, while the positions are random. If $m_a = 2$, we call the swap operator local exploitation, and if $m_a \geq 3$, the swap operator is regarded as global exploration. Fig. \ref{dvsswap} gives the function of the swap transformation graphically when $m_a = 2$.
\\
(2) Shift transformation
\begin{eqnarray}
x_{k+1}= A^{shift}_{k}(m_b) x_{k},
\end{eqnarray}
where, $A^{shift}_{k} \in \mathcal{Z}^{n \times n}$ is called shift permutation matrix, $m_b$ is a constant integer called shift factor to control the maximum length of consecutive positions to be shifted. By the way, the selected position to be shifted after and positions to be shifted are chosen randomly. Similarly, shift transformation is called local exploitation and global exploration when $m_b = 1$ and $m_b \geq2 $  respectively. To make it more clearly, if $m_b = 1$, we set position 2 to be shifted after position 6, as described in Fig. \ref{dvsshift}.
\\
(3) Symmetry transformation
\begin{eqnarray}
x_{k+1}= A^{sym}_{k}(m_c) x_{k},
\end{eqnarray}
where, $A^{sym}_{k} \in \mathcal{Z}^{n \times n}$ is called symmetry permutation matrix, $m_c$ is a constant integer called symmetry factor to control the maximum length of subsequent positions as center. By the way, the component before the subsequent positions and consecutive positions to be symmetrized are both created randomly. Considering that the symmetry transformation can make big change to current solution, it is intrinsically called global exploration. For instance, if $m_c = 0$, let choose the position 3, then the subsequent position or the center is \{$\emptyset$\}, the consecutive positions $\{4,5\}$ with components $(3,1)$, and the function of symmetry transformation is given in Fig. \ref{dvssymmetry}.
\\
(4) Substitute transformation
\begin{eqnarray}
x_{k+1}= A^{sub}_{k}(m_d) x_{k},
\end{eqnarray}
where, $A^{sub}_{k} \in \mathcal{Z}^{n \times n}$ is called substitute permutation matrix, $m_d$ is a constant integer called substitute factor to control the maximum number of positions to be substituted. By the way, the positions are randomly created. If $m_d = 1$, we call the substitute operator local exploitation, and if $m_d \geq 2$, the substitute operator is regarded as global exploration. Fig. \ref{dvssubstitute} gives the function of the substitute transformation vividly when $m_d = 1$.
\subsection{Theoretical analysis}
We give the definition of a global minimum for integer optimization as follows
\begin{eqnarray}\label{eq13}
f( x^{*})  \leq f(x), \exists x^{*} \in \mathcal{X}, \forall\;  x \in \mathcal{X},
\end{eqnarray}
where, $\mathcal{X} \subset \mathcal{Z}^n$ is the feasible space. 
If (\ref{eq13}) is satisfied, we say that $x^{*}$ is a global minimizer.

To give the convergence performance of the proposed discrete STA, we should introduce the general form of random search methods described by  
\begin{algorithmic}[1]
\State Select a starting point $x_{0} \in \mathcal{X}$, and set $k \gets 0$
\Repeat
   \State Generate a candidate solution $ x^{\prime}_{k} \in N(x_{k})\bigcap \mathcal{X}$
   \If{$f(x^{\prime}_{k})< f(x_{k})$}
   \State $x_{k+1} \gets x^{\prime}_{k}$
   \Else
   \State {$x_{k+1} \gets x_{k}$}
   \EndIf
   \State $k \gets k +1 $
\Until{the specified termination criterion is met}
\end{algorithmic}
where, the set $N(x_{k})$ consists of the neighbors of the point $x_{k}$.

It was proved that the above general random search methods can converge in probability to an optimal solution, and it was also demonstrated that, to guarantee the convergence of
an optimization algorithm, the criterion to accept a new solution is a key point \cite{andradottir1996,andradottir1999,baba}. The ``greedy criterion" (always accept a better new solution) is sufficient to guarantee the convergence; nevertheless, accepting a relative worse solution as suggested in simulated annealing can also achieve asymptotic convergence \cite{van}.
In particular, random search methods differ in the choice of the neighborhood structure $N(x_{k})$, and they influence the rate of convergence.

\begin{theorem}
The sequence generated by discrete STA can converge to a global minimizer in probability.
\end{theorem}
\begin{IEEEproof}
On the one hand, the discrete STA is a special case of the random search methods, since the ``greedy criterion" is used as an external archive to keep the incumbent best solution, which can guarantee convergence of the proposed algorithms. 

On the other hand, we have to show that the algorithm can capture a global minimizer probabilistically. Let suppose $x^{*} = (a_1, \cdots, a_n)$ is a global minimum solution, and $ x_N = (b_1, \cdots, b_n)$ is the \textit{N}th best solution. If $x_N =  x^{*}$, according to the ``greedy criterion",
$f( x_k) = f(x_N), \forall ~ k > N$, which means that it converges to $x^{*}$. Otherwise, there must exist a transformation, either swap, shift, symmetry or substitute, such that $x_{N+l_1} = (a_1, \cdots, b_n)$, which means that after $l_1$ iterations, $b_1$ will be changed into $a_1$. If $f(x_{N+l_1}) < f(x_N)$, the $x_{N+l_1}$ is kept as incumbent best for next iteration, else, 
$x_{N+l_1}$ is also kept as incumbent best for next iteration in probability. Following the similar way, after at most $l_2 + \cdots + l_n$ iterations, $x_{N+l_1 + \cdots + l_n}$ will be $(a_1, \cdots, a_n)$ in probability.
\end{IEEEproof}

\subsection{Parameter selection}
In the state transformations, there are four factors to control the intensity between local search and global search. For simplicity and efficiency, the swap, shift and substitute operators are
taken as local search, and the symmetry operator is considered as global search; therefore, we consistently make $m_a = 2, m_b = 1, m_d = 1$ and $m_c = 0$.

On the other hand, the restore probability $p_1$ and the risk probability $p_2$ play a significant role in the discrete STA, as described by the above theorems.
To view the importance of the parameters, we arrange a Monte Carlo simulation study.

Considering the following optimization problem
\begin{eqnarray}
&&\min f^{*}
\end{eqnarray}
where, $f^{*}$ is created by
\begin{algorithmic}[1]
\State {Initialize $f^{*} \gets 0.5, f \gets f^{*}$}
\Repeat
    \If{ $f<$ $r_1$}
    \State {$f$ $\gets$ $r_1$}
    \ElsIf {$r_2 < p_2$} \Comment{risk in probability}
    \State {$f$ $\gets$ $r_1$}
    \EndIf
    \If{ $f^{*}< f$ } \Comment{greedy criterion}
    \State {$f^{*}$ $\gets$ $f$}
    \EndIf
    \If{$r_3 < p_1$} \Comment{restore in probability}
    \State {$f$ $\gets$ $f^{*}$}
    \EndIf
\Until{the maximum number of iterations is met}
\end{algorithmic}
here, $r_1,r_2,r_3$ are uniformly random numbers in $(0,1)$.

We test various groups of $(p_1,p_2)$ for the experiment, in which, the maximum number of iterations is $1e3$, and $1e4$ runs are carried out for each group. The experimental results are shown in
Table \ref{montecarlo}. Without loss of generality, the group $(p_1,p_2) = (0.1,0.1)$ is adopted in this paper for the following study due to its good performance and simplicity.
\begin{table*}[!htbp]
\begin{threeparttable}[b]
\renewcommand{\arraystretch}{1.3}
\caption{A Monte Carlo simulation study}
\label{montecarlo}
\centering
\begin{tabular}{{cccccc}}
\hline
($p_1$ $\setminus$ $p_2$) & 0.1 & 0.3 & 0.5 & 0.7 & 0.9 \\
\hline
0.1&\textbf{9.9254e-4$\pm$ 9.8274e-4}\tnote{1}& 0.0010$\pm$9.8255e-4 & 0.0010 $\pm$ 0.0010 & \textbf{9.9027e-4 $\pm$ 9.9498e-4} &  \textbf{9.7423e-4$\pm$9.8163e-4}\\
0.3&0.0010$\pm$9.9359e-4 & 0.0010 $\pm$ 0.0010 & 0.0010 $\pm$ 0.0010& \textbf{9.8736e-4 $\pm$ 9.8159e-4} & 9.9254e-4 $\pm$ 0.0010\\
0.5&0.0010$\pm$9.9755e-4 & \textbf{9.9249e-4 $\pm$9.7547e-4} & 0.0010 $\pm$ 0.0010& 0.0010 $\pm$ 0.0010& 0.0010 $\pm$ 0.0010\\
0.7&0.0010 $\pm$ 0.0010  & 0.0010 $\pm$ 0.0010 & 0.0010 $\pm$ 0.0010 & \textbf{9.9456e-4 $\pm$ 9.8730e-4} & \textbf{9.9609e-4 $\pm$  9.9991e-4}\\
0.9& 0.0010 $\pm$ 9.9637e-4 & 0.0010 $\pm$ 9.9667e-4 & \textbf{9.8258e-4 $\pm$ 9.8261e-4} & 0.0010 $\pm$ 0.0010& 0.0010 $\pm$ 0.0010\\
\hline
\end{tabular}
\begin{tablenotes}
\item [1] indicates mean $\pm$ standard deviation
\end{tablenotes}
\end{threeparttable}
\end{table*}
\section{Implementation of the discrete STA}
The above discrete STA are essentially for unconstrained discrete optimization problem. To realize the optimal design of water distribution networks, we have to deal with some constraints.
For the equality constraints on continuity of flow and conservation of energy, there exist some hydraulic analysis software packages such as EPANET \cite{rossman1994}, KYPIPE \cite{wood1980}, in which the continuity and energy constraints are automatically satisfied.
Considering that the continuity equations are linear, we can first fix some of pipe flows as known to solve the linear equations and then substitute them into the energy equations, which can reduce the computational complexity of solving continuity equations (linear) and energy equations (nonlinear) simultaneously. It is not difficult to imagine that the number of nonlinear equations equals to that of the simple closed loops or required independent paths in a network, and then a Newton-Raphson method is used to solve the nonlinear equations. \\
\indent For the minimum pressure head constraints, the most commonly used technique is the penalty function method, adding a penalty term when the corresponding constraint is violated. For example, the following scheme
\begin{eqnarray}
f_{penal} = pc \sum_{i=1}^{NP}\max\{0, H_{i\min} - H_i\}^{\rho}
\end{eqnarray}
where, $pc$ is the penalty coefficient, and ${\rho}$ is normally 1 or 2 (${\rho} = 1$ in this study). Finally, the total cost is
\begin{eqnarray}
f_{cost} = f_{obj} +  f_{penal}.
\end{eqnarray}
\indent A brief description of the steps using discrete STA is given in the following
\begin{enumerate}
  \item Creat initial $Best$ solution. Generate a group of candidate solutions randomly (the size is the search enforcement, \textit{SE}) and then select the fittest solution. Let $Best^{*} = Best$ and store $Best^{*}$.
  \item Update the $Best$. Use swap transformation to generate a group of candidate solutions on the basis of $Best$. If the fittest of the candidate solutions is better than $Best$, then accept the fittest solution as $Best$; otherwise, accept the fittest solution as $Best$ in a probability $p_2$. Similar procedures are adaptive to shift, symmetry and substitute transformations.
  \item Update the $Best^{*}$. The $Best^{*}$ is updating only when $Best$ is better than $Best^{*}$.
  \item Restore the $Best$. The $Best$ is restored to $Best^{*}$ in a probability $p_1$.
  \item Go back to repeat step 2 until the stopping criterion is met.
\end{enumerate}
\begin{remark}
We should notice that once a solution is given, then the flow in each pipe is determined by solving the nonlinear equations, and then we can evaluate whether the minimum pressure head is satisfied and decide the corresponding penalty term to each head pressure constraint.
\end{remark}
\section{Case studies}
We investigate the performance of the proposed discrete STA by three well-known water distribution networks, namely, the Two-Loop network, the Hanoi network and the New York network.
We first give a detailed study of the Two-Loop problem to show that the network system with eight unknowns governed by six linear equations and two nonlinear equations can be reduced to only two unknowns governed by two nonlinear equations. This Two-Loop case is also fully trained to study the effect of penalty coefficient and search enforcement on the performance of the algorithm.
Based on the experience gained from the case, the known best solutions for the other two networks are also achieved by the algorithm.

%For the two-loop network, the penalty function method with different schemes and another strategy based on the objective function and the constraint violation are applied respectively.
\subsection{two-loop network}
\begin{figure}[!htbp]
  \centering
  % Requires \usepackage{graphicx}
  \includegraphics[width=6cm]{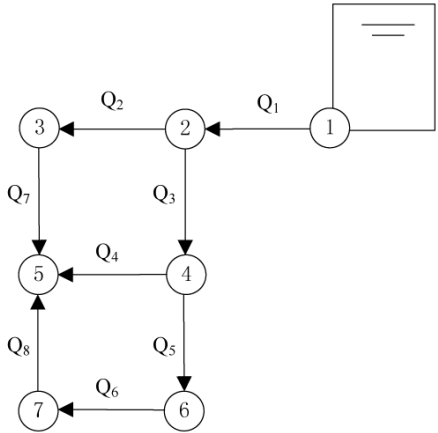}\\
  \caption{two loop network}\label{twoloop}
\end{figure}

The layout of the two-loop network is given in Fig. \ref{twoloop}. There are a single reservoir at a $210$-m fixed head and eight
pipes all with $1000$-m long. The node data and cost data are given in Table \ref{nodedatatwoloop} and Table \ref{costdatatwoloop}, and the minimum acceptable pressure requirements are all 30-m above the ground level. The Hazen-Williams coefficient $C$ is assumed to be $130$ for the two-loop network.

\begin{table}[!htbp]
\renewcommand{\arraystretch}{1.3}
\caption{Node data for the two loop network}
\label{nodedatatwoloop}
\centering
\begin{tabular}{{ccc}}
\hline
Node  & Demand($m^3/h$)  & Ground level($m$)\\
\hline
1 & -1120.0 & 210.00\\
2 & 100.0   & 150.00\\
3 & 100.0   & 160.00\\
4 & 120.0   & 155.00\\
5 & 270.0   & 150.00\\
6 & 330.0   & 165.00\\
7 & 200.0   & 160.00\\
\hline
\end{tabular}
\end{table}
\begin{table}[!htbp]
\begin{threeparttable}[b]
\renewcommand{\arraystretch}{1.3}
\caption{Cost data for the two loop network}
\label{costdatatwoloop}
\centering
\begin{tabular}{{cccccc}}
\hline
No. & Diameter (in.)\tnote{2}   & Cost (\$/m) & No. & Diameter (in.)   & Cost(\$/m)\\
\hline
    1&    1&       2&      8&   12&     50\\
    2&    2&      5&     9&   14&     60\\
    3&    3&     8&    10&   16&    90\\
    4&    4&    11&    11&   18&   130\\
    5&    6&    16&    12&   20&     170\\
    6&    8&   23&    13&   22&   300\\
    7&   10&    32&   14&   24&   550\\
\hline
\end{tabular}
\begin{tablenotes}
\item [2] 1 in. = 2.54 cm
\end{tablenotes}
\end{threeparttable}
\end{table}
In this case study, we give an illustrative procedure of how to reduce the complexity of solving the linear and nonlinear equations.
The flow continuity equations of the two-loop network are given as follows
\begin{equation}
\left\{ \begin{aligned}
& -Q_1 + Q_2 + Q_3 + DM_2 = 0 \\
& -Q_2 + Q_7 + DM_3 = 0\\
& -Q_3 + Q_4 + Q_5 + DM_4 = 0\\
& -Q_7 - Q_8 - Q_4 + DM_5 = 0\\
& -Q_5 + Q_6 + DM_6 = 0\\
& -Q_6 + Q_8 + DM_7 = 0
\end{aligned} \right.
\end{equation}
Let $Q_4, Q_6$ be fixed, then
\begin{equation}
\left\{ \begin{aligned}
& Q_1 = DM_2 + DM_3 + DM_4 + DM_5 + DM_6 + DM_7\\
& Q_2 = DM_3 + DM_5 + DM_7 - Q_4 - Q_6\\
& Q_3 = DM_4 + DM_6 + Q_4 + Q_6\\
& Q_5 = DM_6 + Q_6 \\
& Q_7 = DM_5 + DM_7 - Q_4 - Q_6\\
& Q_8 = Q_6 - DM_7
\end{aligned} \right.
\end{equation}
The energy conservation equations can be formulated as
\begin{small}
\begin{equation} \label{eq:19}
\left\{ \begin{aligned}
&r_3 Q_3|Q_3|^{\alpha-1} \!+ \!r_4 Q_4|Q_4|^{\alpha-1} \!- \!r_7 Q_7|Q_7|^{\alpha-1}\!-\! r_2 Q_2|Q_2|^{\alpha-1}\!=\!0\\
&r_5 Q_5|Q_5|^{\alpha-1} \!+\! r_6 Q_6|Q_6|^{\alpha-1} \!+\! r_8 Q_8|Q_8|^{\alpha-1} \!-\! r_4 Q_4|Q_4|^{\alpha-1}\!=\!0
\end{aligned} \right.
\end{equation}
\end{small}
and the head loss equations
\begin{equation}
\left\{ \begin{aligned}
& H_2 = Head - r_1 Q_1|Q_1|^{\alpha-1} - G_2 \geq H_{2\min}\\
& H_3 = H_2 - r_2 Q_2|Q_2|^{\alpha-1} - G_3 \geq H_{3\min}\\
& H_4 = H_2 - r_3 Q_3|Q_3|^{\alpha-1} - G_4 \geq H_{4\min}\\
& H_5 = H_4 - r_4 Q_4|Q_4|^{\alpha-1} - G_5 \geq H_{5\min}\\
& H_6 = H_4 - r_5 Q_5|Q_5|^{\alpha-1} - G_6 \geq H_{6\min}\\
& H_7 = H_6 - r_6 Q_6|Q_6|^{\alpha-1} - G_7 \geq H_{7\min}\\
\end{aligned} \right.
\end{equation}
where, $G_i(i = 2, \cdots,7)$ is the ground level.

\begin{remark}
It should be noted that we only need to solve the nonlinear system (\ref{eq:19}) with two unknowns ($Q_4, Q_6$).
\end{remark}

For the Two-Loop network, we have to select a diameter to each pipe, and for each pipe there are 14 choices. It is not difficult to imagine that when choosing a numerical order (No.), it corresponds to an exact diameter. That is the reason why the discrete STA use the permutation of $\{1,2,\cdots,n\}$ as its decision variables and all the intelligent operators are operated on a certain permutation.

Next, we conduct a eremitical study of the the two loop network by the proposed discrete STA to investigate the influence of the remained parameters, namely, the search enforcement ($SE$) and
the penalty coefficient ($pc$). We set $SE$ to be $0.5,1,2,3$ and $4$ times of the dimension of decision variable. Considering that the average cost times the average pipe length is 1.0335e5 and the average of the minimum pressure heads is 30, the order of magnitude for $pc$ is set at 1e4. On this situation, $pc$ is fixed at $1e4,2e4,4e4,8e4$ and $1e5$, or increases from $1e4$ to $1e5$ in a linear way. The maximum number of iterations is set at $2e2$, and a total of $20$ runs are executed for each group of search enforcement $SE$ and penalty coefficient $pc$.
\begin{table*}[!htbp]
\begin{threeparttable}[b]
\renewcommand{\arraystretch}{1.3}
\caption{A empirical study of the two loop network}
\scriptsize
\label{empiricalstudy}
\centering
\begin{tabular}{{p{1.1cm}p{2.2cm}p{2.2cm}p{2.2cm}p{2.2cm}p{2.2cm}p{2.2cm}}}
\hline
($SE$ $\setminus$ $pc$) & 1e4 & 2e4 & 4e4 & 8e4 & 1e5 & 1e4  $\rightarrow$ 1e5 \\
\hline
4& 4.2978e5 $\pm$ 1.4882e4 (55\%)\tnote{3} & 4.3631e5 $\pm$ 1.3394e4 (85\%) & 4.5184e5 $\pm$ 2.3575e4 (95\%) & 4.5063e5 $\pm$ 1.7400e4 (95\%) & 4.4190e5 $\pm$ 1.6121e4 (95\%) & 4.4368e5 $\pm$ 1.8563e4 (90\%)  \\
8& 4.2195e5 $\pm$ 1.4853e4 (65\%) & \textbf{4.3181e5 $\pm$ 1.3870e4 (85\%)} & 4.3526e5 $\pm$ 1.2721e4 (90\%)  & 4.3577e5 $\pm$ 1.2903e4 (95\%) & 4.4085e5 $\pm$ 1.5853e4 (90\%) & 4.3620e5 $\pm$ 1.5702e4 (80\%)\\
16& 4.2682e5 $\pm$ 1.2946e4 (75\%)& 4.3340e5 $\pm$ 1.5347e4 (80\%) & 4.3410e5 $\pm$ 1.2004e4 (90\%)  & 431550 $\pm$ 1.4406e4 (100\%) & 4.3458e5 $\pm$ 1.4992e4 (90\%) & 433600 $\pm$ 1.4207e4 (100\%)\\
24& 4.2380e5 $\pm$ 1.2756e4 (75\%)& 4.3193e5 $\pm$ 1.2898e4 (95\%) & 4.3555e5 $\pm$ 1.5049e4 (90\%)  & 440950 $\pm$ 1.4417e4 (100\%) & 432600 $\pm$ 1.5398e4 (100\%) & 4.3073e5 $\pm$ 1.3252e4 (95\%)\\
32& 4.2686e5 $\pm$ 1.5549e4 (55\%)& 4.3046e5 $\pm$ 1.5523e4 (80\%) & 4.3376e5 $\pm$ 1.3900e4 (95\%)  & 451400 $\pm$ 5.4648e4 (100\%) & 4.3600e5 $\pm$ 1.7731e4 (85\%) & 4.3305e5 $\pm$ 1.4657e4 (95\%)\\
\hline
\end{tabular}
\begin{tablenotes}
\item [3] indicates the percentage of feasible solutions
\end{tablenotes}
\end{threeparttable}
\end{table*}

\begin{figure*}[!htbp]
  \centering
  % Requires \usepackage{graphicx}
  \includegraphics[width=8cm]{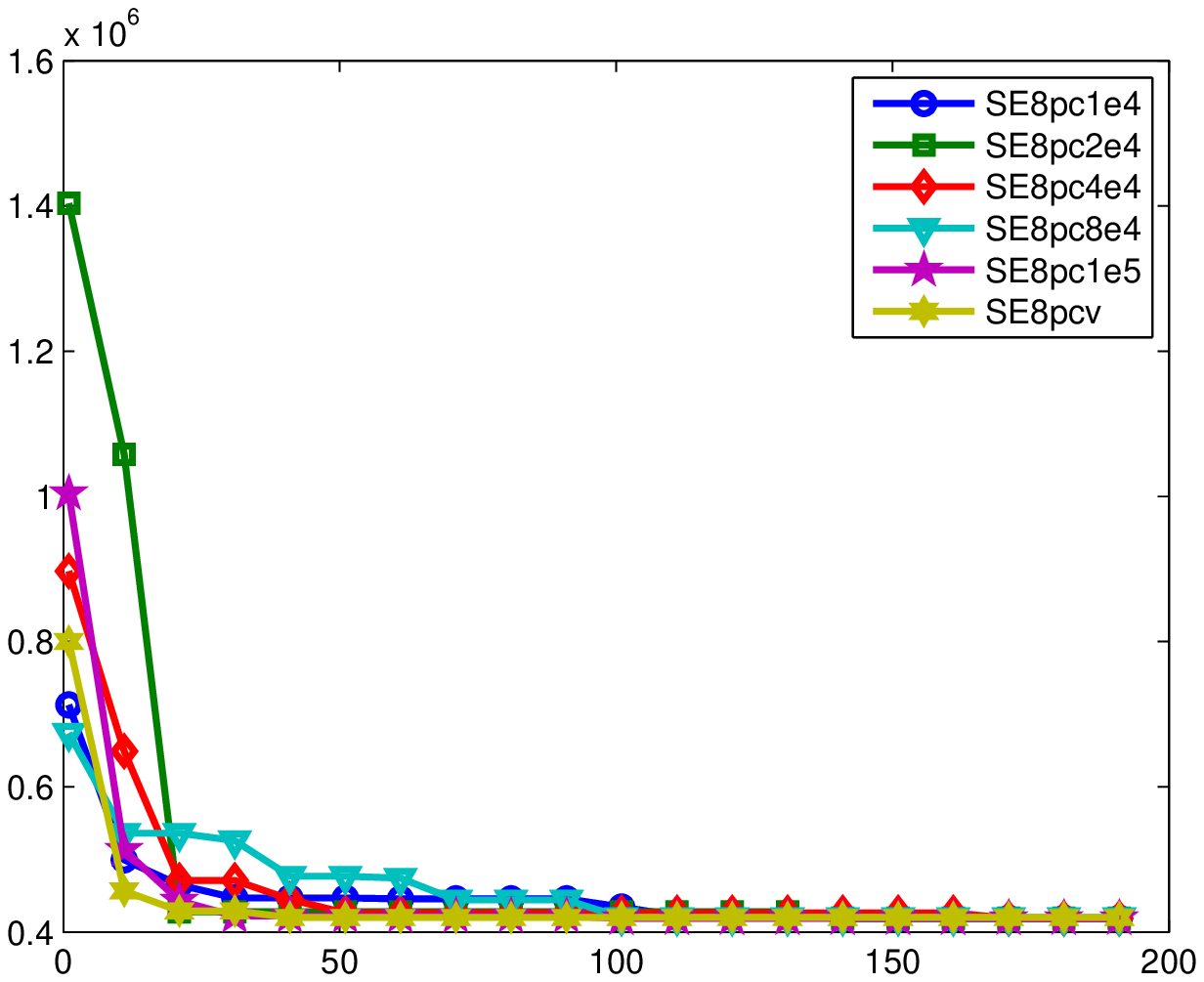}
  \includegraphics[width=8cm]{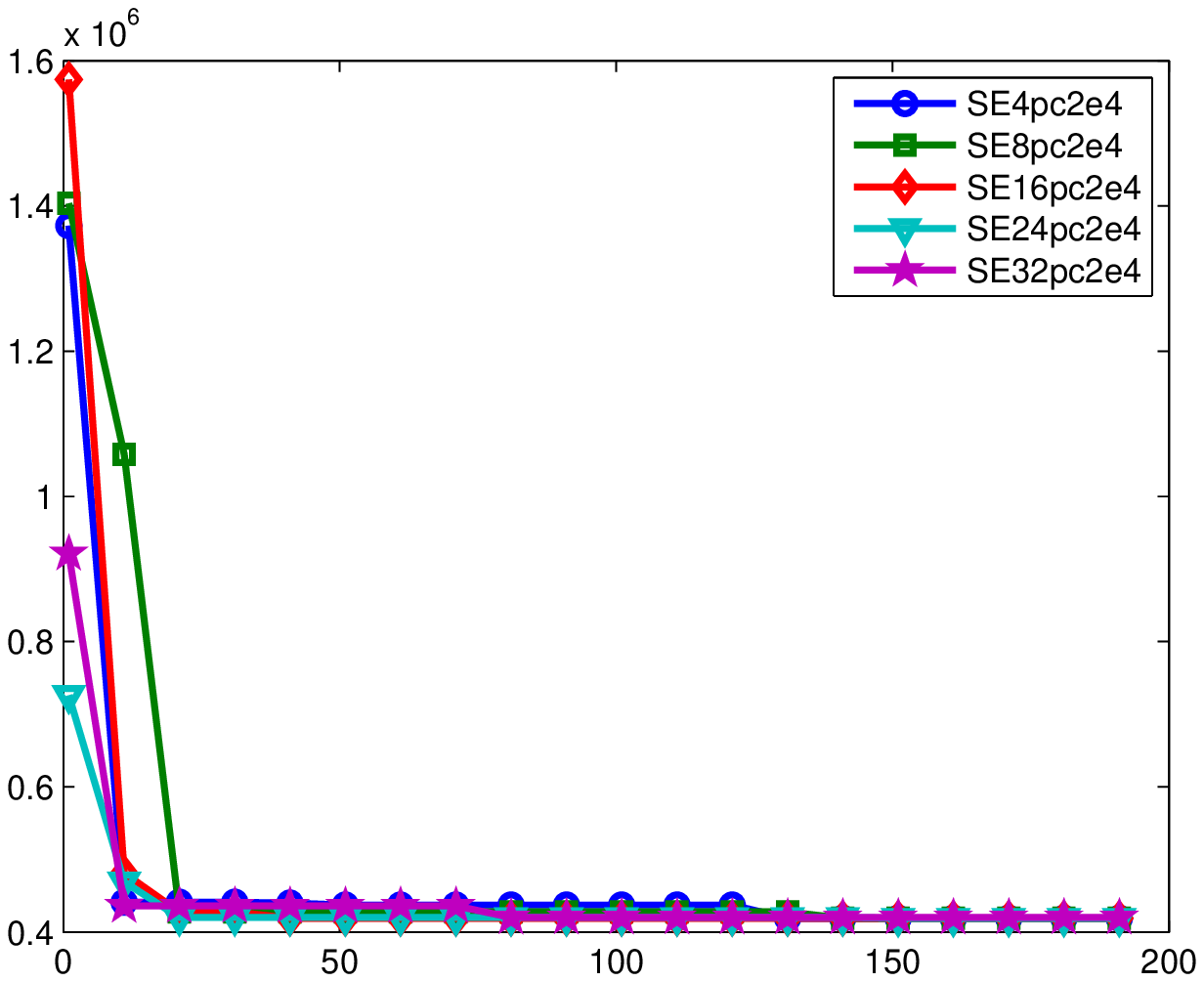}
  \caption{Iterative curves of best solutions when $SE = 8$ and $pc = 2e4$ for the Two-Loop problem, respectively}
  \label{se8pc2e4}
\end{figure*}

As can be seen from Table \ref{empiricalstudy}, for a fixed $SE$, the search ability is declining as the $pc$ increases, but the feasibility rate increases simultaneously with the $pc$. For a
fixed $pc$, the search ability is increasing as the $SE$ increases from $4$ to $8$ but declining as the $SE$ increases any more. When the $pc$ varies in the iterative process, the performance is not the best but much more satisfactory than a constant one to some extent. By observation, we can find that setting $SE$ to be the dimensionality of the decision variable is a good choice, and in this setting environment, $pc = 2e4$ is a good penalty coefficient. Fig. \ref{se8pc2e4} gives the iterative curves of the gained best solutions when $SE = 8$ and $pc = 2e4$, respectively.
It should be emphasized that best solutions are all $419,000$.
\begin{remark}
Under the circumstance, the minimum function evaluations to achieve the best known solution is $2048$, which takes up $0.0001387\%$ of all possible combinations ($14^8 = 1.4758e9$).
\end{remark}

Table \ref{twoloopsolution} gives the best solutions gained by various algorithms, and it can be found that STA can achieve the best known solution in this case. It should be noted that the same solution was also achieved by GA \cite{savic1997}, SA \cite{cunha1999} and HS \cite{geem2006} with function evaluations at $250,000$, $70,000$ and $5,000$, respectively. Although the solution in \cite{kessler1989} is even better, it should be noted that it brings pipe segments. The pressure heads for the Two-Loop network obtained by various algorithms are given in Table \ref{twolooppressureheads}.
\begin{table*}[!htbp]
\caption{Solutions for the two loop network}
\label{twoloopsolution}
\renewcommand{\arraystretch}{1.3}
\centering
\footnotesize
\begin{tabular}{{cccccc}}
\hline
Pipe  & Alperovits and Shamir \cite{alperovits1977}   & Goulter et al. \cite{goulter1986} & Kessler and Shamir \cite{kessler1989} & STA (fixed)  & STA (variable)\\
\hline
                                             \\
    1&    20&       20&      18&   18&     18\\
     &    18&       18&                    \\
    2&     8&      10&     12&    10&     10\\
     &     6&        &     10&      &     \\
    3&    18&      16&     16&    16&     16\\
    4&     8&       6&     3&      4&     4\\
     &     6&       4&     2&       &      \\
    5&    16&      16&    16&     16&     16\\
     &      &      14&    14&       &       \\
    6&    12&      12&    12&     10&     10\\
     &    10&      10&    10&       &      \\
    7&     6&      10&    10&     10&     10\\
     &      &       8&     8&       &     \\
    8&     6&       2&     3&      1&     1\\
     &     4&       1&     2&       &     \\
Cost(\$) & 497,525 & 435,015 & 417,500 & 419,000 & 419,000 \\
\hline
\end{tabular}
\end{table*}

\begin{table*}[!htbp]
\caption{Pressure Heads for the Two-Loop network}
\label{twolooppressureheads}
\renewcommand{\arraystretch}{1.3}
\centering
\footnotesize
\begin{tabular}{{cccccc}}
\hline
Node  & Alperovits and Shamir \cite{alperovits1977}   & Goulter et al. \cite{goulter1986} & Kessler and Shamir \cite{kessler1989} & STA (fixed and variable)\\
\hline
2   & 53.96 & 54.30 & 53.26 & 53.24 & \\
3   & 32.32 & 33.19 & 30.08 & 30.49 & \\
4   & 44.97 & 44.19 & 43.64 & 43.44 & \\
5   & 32.31 & 32.32 & 30.10 & 33.78 & \\
6   & 31.19 & 31.19 & 30.08 & 30.43 & \\
7   & 31.57 & 31.57 & 30.09 & 30.54 & \\
\hline
\end{tabular}
\end{table*}

\subsection{Hanoi network}
\begin{figure*}[!htbp]
  \centering
  % Requires \usepackage{graphicx}
  \includegraphics[width=14cm]{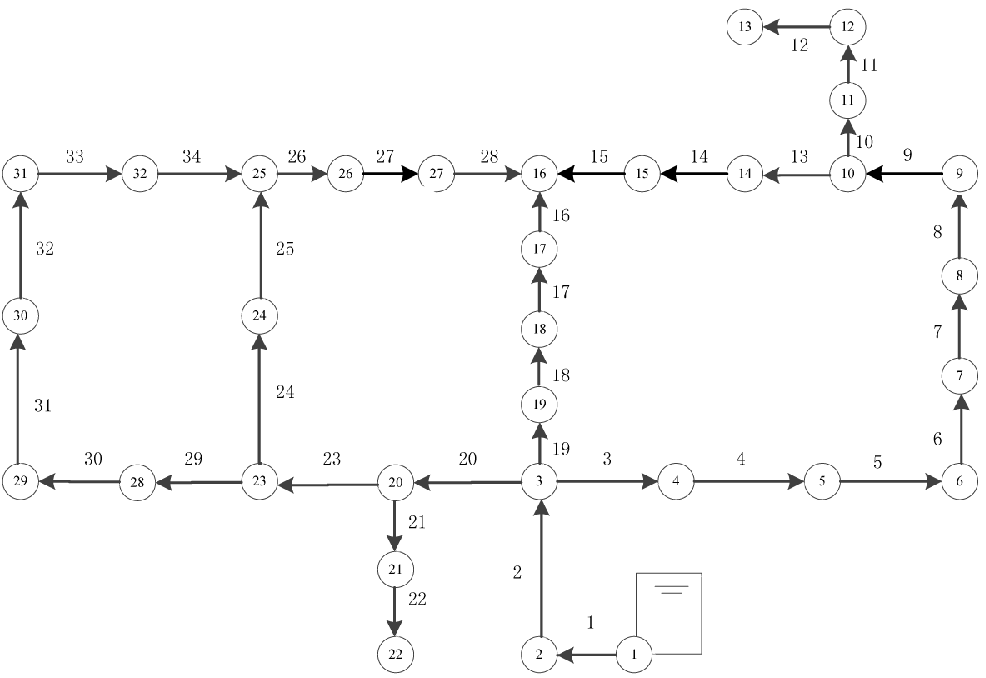}\\
  \caption{The Hanoi network}\label{hanoi}
\end{figure*}

The layout of the Hanoi network is given in Fig. \ref{hanoi}. There are $32$ nodes, $34$ pipes and $3$ loops in this network system. At node 1, there exists a reservoir with
a 100-m fixed head. The cost data, and pipe and node data are given in Table \ref{costdatahanoi} and Table \ref{pipenodedatahanoi}, respectively. The minimum acceptable pressure requirements at all nodes are also fixed at 30 m and the Hazen-Williams coefficient $C$ is assumed to be $130$ as well.
\begin{table}[!htbp]
\renewcommand{\arraystretch}{1.3}
\caption{Cost data for the Hanoi network}
\label{costdatahanoi}
\centering
\footnotesize
\begin{tabular}{{ccc}}
\hline
No. & Diameter (in.)   & Cost (\$/m)\\
\hline
1 & 12 & 45.726\\
2 & 16 & 70.400\\
3 & 20 & 98.387\\
4 & 24 & 129.333\\
5 & 30 & 180.748\\
6 & 40 & 278.280\\
\hline
\end{tabular}
\end{table}
\begin{table*}[!htbp]
\renewcommand{\arraystretch}{1.3}
\caption{Pipe and node data for the Hanoi network}
\label{pipenodedatahanoi}
\centering
\footnotesize
\begin{tabular}{{cccccccc}}
\hline
Pipe  & Length ($m$)  & Pipe & Length ($m$) & Node  & Demand ($m^3/h$) & Node  & Demand ($m^3/h$) \\
\hline
1 & 100 & 18& 800 & 1 & -19940 & 18& 1345\\
2 & 1350& 19& 400 & 2 & 890    & 19& 60 \\
3 & 900 & 20& 2200& 3 & 850    & 20& 1275\\
4 & 1150& 21& 1500& 4 & 130    & 21& 930\\
5 & 1450& 22& 500 & 5 & 725    & 22& 485\\
6 & 450 & 23& 2650& 6 & 1005   & 23& 1045\\
7 & 850 & 24& 1230& 7 & 1350   & 24& 820\\
8 & 850 & 25& 1300& 8 & 550    & 25& 170\\
9 & 800 & 26& 850 & 9 & 525    & 26& 900\\
10& 950 & 27& 300 & 10& 525    & 27& 370\\
11& 1200& 28& 750 & 11& 500    & 28& 290\\
12& 3500& 29&1500 & 12& 560    & 29&360\\
13& 800 & 30&2000 & 13& 940    & 30&360\\
14& 500 & 31&1600 & 14& 615    & 31&105\\
15& 550 & 32&150  & 15& 280    & 32&805\\
16& 2730& 33&860  & 16& 310    & - & -\\
17& 1750& 34&950  & 17& 865    & - & -\\
\hline
\end{tabular}
\end{table*}

From the experience gained from the training of the Two-Loop network, the search enforcement $SE$ does not affect the performance of the discrete STA explicitly, but the penalty coefficient $pc$ plays a significant role in the search ability and the solution feasibility, and a good penalty coefficient can be evaluated from the order of magnitude the same as the average pipe length times the minimum pressure heads.

\begin{figure*}[!htbp]
  \centering
  % Requires \usepackage{graphicx}
  \includegraphics[width=8cm]{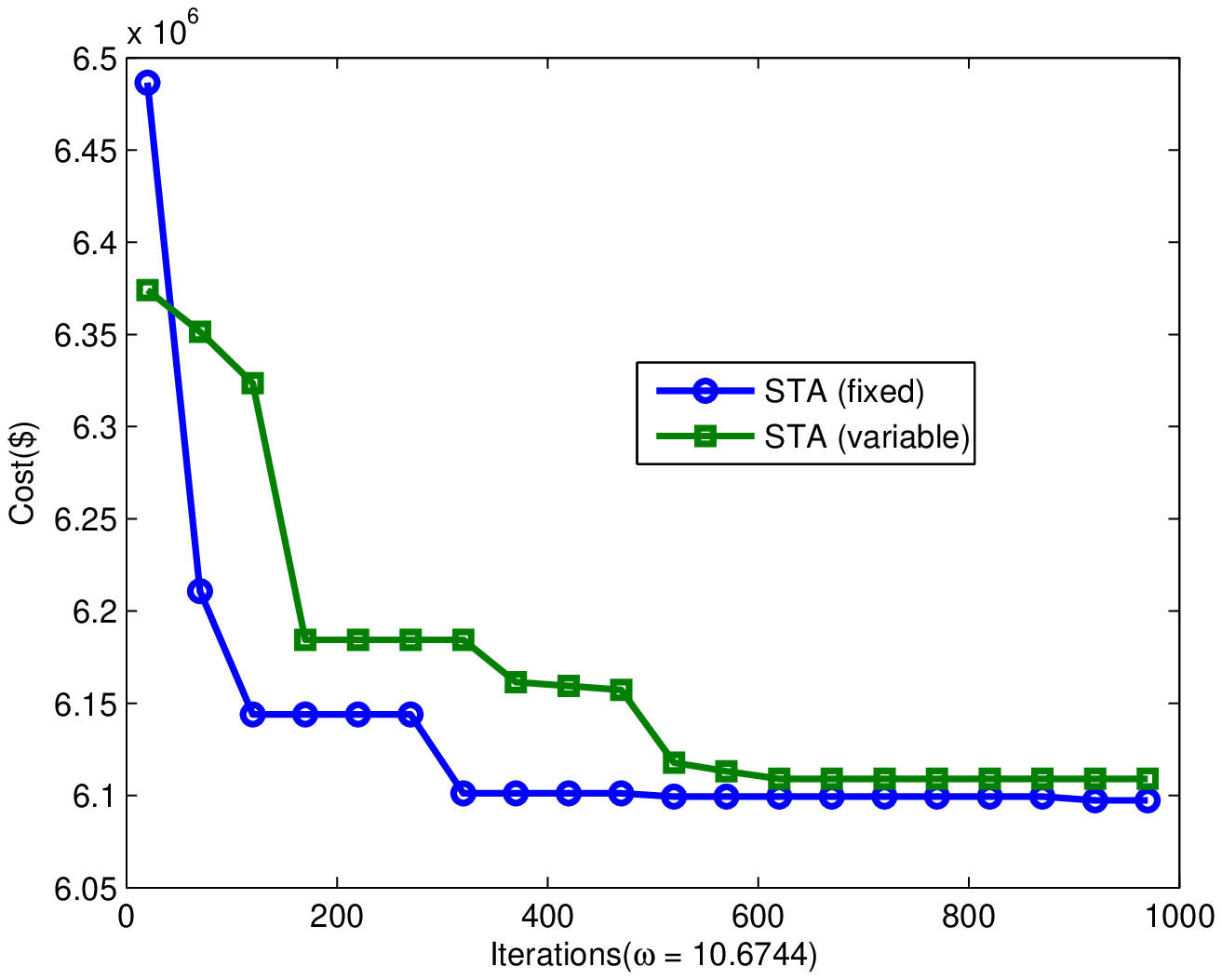}
  \includegraphics[width=8cm]{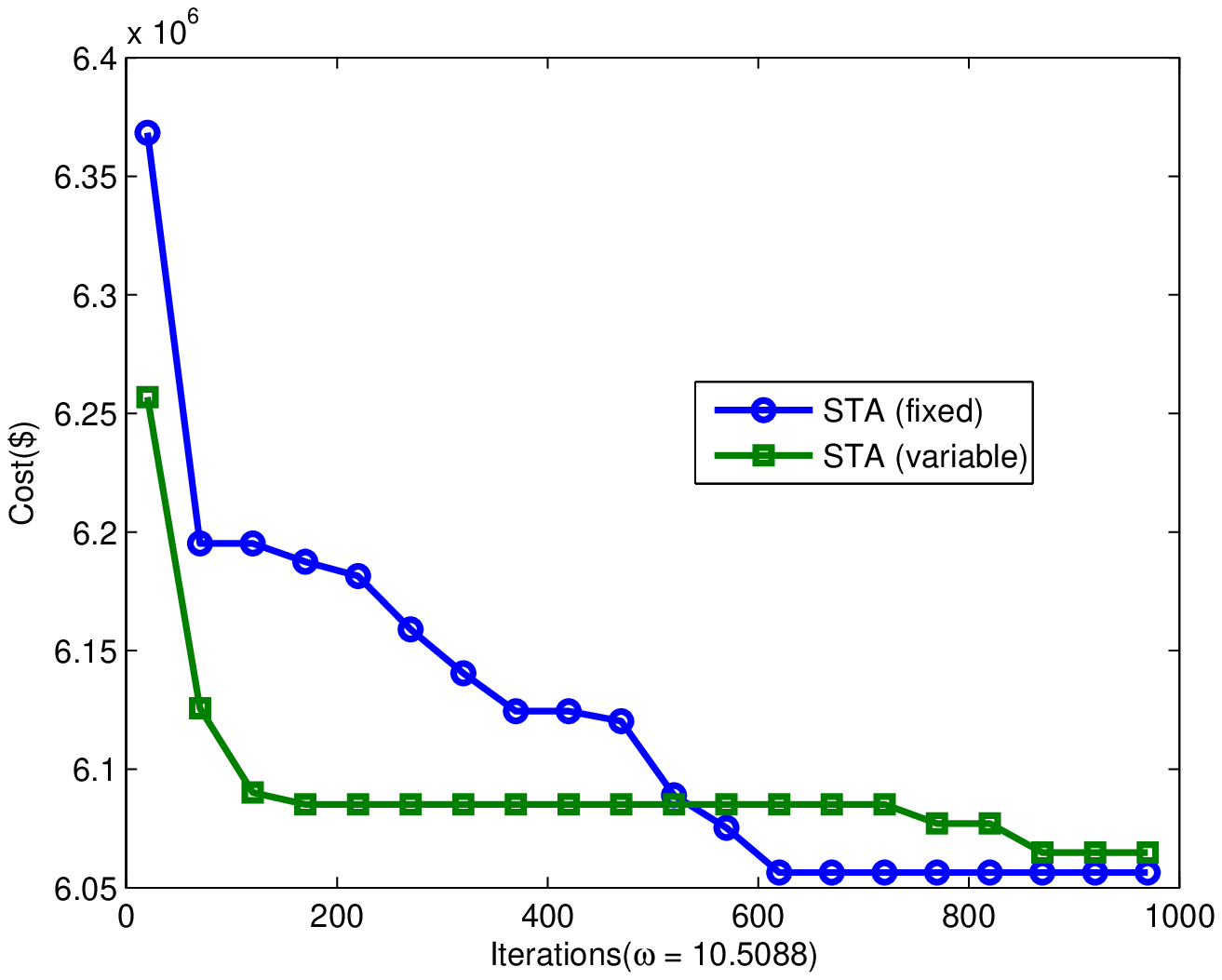}\\
  \caption{Iterative curves of best solutions using STA for the Hanoi problem when $\omega$ is 10.6744 and 10.5088, respectively}\label{hanoicurveruns}
\end{figure*}

For the Hanoi network, the search enforcement $SE$ is set at $20$, and the penalty coefficient $pc$ is fixed at $4e4$, or varies from $4e4$ to $1e5$ in a linearly increasing way.
The maximum number of iterations is set at $1e3$, and a total of $20$ runs are executed for both fixed and variable $pc$. Fig. \ref{hanoicurveruns} gives the
iterative curves of best solutions and changes in $20$ runs for the Hanoi problem with fixed and variable $pc$ respectively. It is shown that the best solution is hit only once by the STA with fixed $pc$.

\begin{remark}
Under the circumstance, the minimum function evaluations to achieve the best known solution is $23,240$, which takes up $8.1114e$-$21\%$ of all possible combinations ($6^{34} = 2.8651e26$).
\end{remark}

Table \ref{hanoisolution} gives the best solutions gained by various algorithms, and it can be found that STA with fixed $pc$ can achieve the best known solution in this case  at the cost of $6.056$ million dollars, while the solution of STA with variable $pc$ get a solution at the cost of $6.065$ million dollars. Savic and Walters \cite{savic1997} used the GA to obtain the solution with 1,000,000 function evaluations. The solution gained by Zecchin et al. \cite{zecchin2006} using ACO need 100,000 function evaluations. The exactly same solution was achieved by SA \cite{cunha1999} and HS \cite{geem2006} as well, with the function evaluations at $53,000$ and $200,000$, respectively. The pressure heads for the Hanoi network obtained by various algorithms are given in Table \ref{hanoipressureheads}.

\begin{table*}[!htbp]
\caption{Solutions for the Hanoi network}
\label{hanoisolution}
\renewcommand{\arraystretch}{1.3}
\centering
\footnotesize
\begin{tabular}{{p{2.0cm}p{1.3cm}p{1.2cm}p{1.2cm}p{1.8cm}p{1.8cm}p{1.8cm}p{1.8cm}}}
\hline
Pipe  & Savic and Walters \cite{savic1997}   &  Zecchin ~et al. \cite{zecchin2006}  & Haghighi ~et al. \cite{haghighi2011}& STA (fixed) $\omega = 10.6744$& \;$\omega = 10.5088$ & STA (variable) $\omega = 10.6744$ &\; $\omega = 10.5088$ \\
\hline
1    & 40 & 40 & 40 & 40 &40& 40& 40\\
2    & 40 & 40 & 40 & 40 &40& 40& 40\\
3    & 40 & 40 & 40 & 40 &40& 40& 40\\
4    & 40 & 40 & 40 & 40 &40& 40& 40\\
5    & 40 & 40 & 40 & 40 &40& 40& 40\\
6    & 40 & 40 & 40 & 40 &40& 40& 40\\
7    & 40 & 40 & 40 & 40 &40& 40& 40\\
8    & 40 & 40 & 40 & 40 &40& 40& 40\\
9    & 40 & 40 & 30 & 40 &40& 30& 30\\
10   & 30 & 30 & 30 & 30 &30& 30& 30\\
11   & 24 & 24 & 30 & 24 &24& 30& 30\\
12   & 24 & 24 & 24 & 24 &24& 24& 24\\
13   & 20 & 20 & 16 & 20 &20& 20& 20\\
14   & 16 & 12 & 12 & 16 &16& 12& 16\\
15   & 12 & 12 & 12 & 12 &12& 12& 12\\
16   & 12 & 12 & 16 & 12 &12& 12& 12\\
17   & 16 & 20 & 20 & 16 &16& 20& 16\\
18   & 20 & 24 & 24 & 24 &20& 20& 24\\
19   & 20 & 20 & 24 & 20 &20& 24& 20\\
20   & 40 & 40 & 40 & 40 &40& 40& 40\\
21   & 20 & 20 & 20 & 20 &20& 20& 20\\
22   & 12 & 12 & 12 & 12 &12& 12& 12\\
23   & 40 & 40 & 40 & 40 &40& 40& 40\\
24   & 30 & 30 & 30 & 30 &30& 30& 30\\
25   & 30 & 30 & 30 & 30 &30& 30& 30\\
26   & 20 & 20 & 20 & 20 &20& 20& 20\\
27   & 12 & 12 & 12 & 12 &12& 12& 12\\
28   & 12 & 12 & 12 & 12 &12& 12& 12\\
29   & 16 & 16 & 16 & 16 &16& 16& 16\\
30   & 16 & 16 & 12 & 16 &12& 16& 12\\
31   & 12 & 12 & 12 & 12 &12& 12& 12\\
32   & 12 & 12 & 16 & 12 &16& 16& 16\\
33   & 16 & 16 & 20 & 16 &16& 16& 16\\
34   & 20 & 20 & 24 & 20 &24& 20& 24\\
Cost(\$ millions) & 6.073 & 6.134 &  6.190 & 6.097 & 6.056 & 6.109 & 6.065 \\
\hline
\end{tabular}
\end{table*}

\begin{table*}[!htbp]
\caption{Pressure Heads for the Hanoi network}
\label{hanoipressureheads}
\renewcommand{\arraystretch}{1.3}
\centering
\footnotesize
\begin{tabular}{{p{1.0cm}p{1.3cm}p{1.2cm}p{1.8cm}p{1.8cm}p{1.8cm}p{1.8cm}}}
\hline
Node  & Savic and Walters \cite{savic1997}   & Haghighi et al. \cite{haghighi2011} & STA (fixed) $\omega = 10.6744$& \;$\omega = 10.5088$  & STA (variable) $\omega = 10.6744$& \;$\omega = 10.5088$ \\
\hline
1 & 100.00 & 100.00 & 100.00 & 100.00  & 100.00 & 100.00\\
2 & 97.16  & 97.08  & 97.14  & 97.17   & 97.14  & 97.17\\
3 & 61.95  & 60.82  & 61.64  & 61.99   & 61.64  & 61.99\\
4 & 57.21  & 56.38  & 56.90  & 57.23   & 57.08  & 57.34\\
5 & 51.33  & 50.88  & 51.02  & 51.31   & 51.42  & 51.56\\
6 & 45.13  & 45.13  & 44.82  & 45.07   & 45.49  & 45.48\\
7 & 43.68  & 43.81  & 43.36  & 43.61   & 44.10  & 44.06\\
8 & 41.93  & 42.28  & 41.63  & 41.85   & 42.48  & 42.37\\
9 & 40.54  & 41.09  & 40.25  & 40.44   & 41.20  & 41.02\\
10& 40.34  & 37.61  & 39.23  & 39.40   & 37.39  & 37.01\\
11& 38.79  & 36.01  & 37.67  & 37.85   & 35.83  & 35.45\\
12& 38.78  & 34.83  & 34.24  & 34.43   & 34.68  & 34.30\\
13& 34.58  & 30.53  & 30.03  & 30.24   & 30.46  & 30.10\\
14& 36.59  & 32.06  & 35.61  & 35.49   & 34.64  & 33.66\\
15& 34.71  & 30.96  & 33.87  & 33.44   & 30.86  & 32.17\\
16& 32.08  & 31.13  & 31.61  & 30.36   & 30.38  & 30.53\\
17& 33.36  & 39.28  & 33.56  & 30.51   & 38.00  & 33.20\\
18& 43.32  & 50.04  & 49.94  & 44.29   & 44.89  & 50.16\\
19& 55.54  & 57.13  & 55.08  & 55.90   & 58.68  & 55.37\\
20& 50.92  & 49.59  & 50.53  & 50.89   & 50.43  & 50.90\\
21& 44.79  & 40.04  & 41.18  & 41.57   & 41.07  & 41.59\\
22& 39.63  & 34.76  & 36.01  & 36.42   & 35.90  & 36.44\\
23& 44.83  & 43.42  & 44.41  & 44.73   & 44.21  & 44.76\\
24& 39.64  & 37.73  & 39.23  & 39.03   & 38.91  & 39.07\\
25& 36.38  & 34.07  & 35.98  & 35.34   & 35.56  & 35.40\\
26& 32.67  & 30.51  & 32.25  & 31.44   & 31.58  & 31.53\\
27& 31.66  & 30.32  & 31.20  & 30.15   & 30.22  & 30.29\\
28& 36.48  & 38.05  & 35.76  & 39.12   & 35.60  & 39.15\\
29& 32.04  & 30.08  & 31.06  & 30.21   & 30.94  & 30.26\\
30& 31.29  & 30.58  & 30.10  & 30.47   & 30.01  & 30.52\\
31& 31.81  & 30.90  & 30.58  & 30.75   & 30.13  & 30.80\\
32& 32.17  & 31.81  & 31.84  & 33.20   & 31.41  & 33.26\\
\hline
\end{tabular}
\end{table*}

\subsection{New York network}
\begin{figure}[!htbp]
  \centering
  % Requires \usepackage{graphicx}
  \includegraphics[width=7cm]{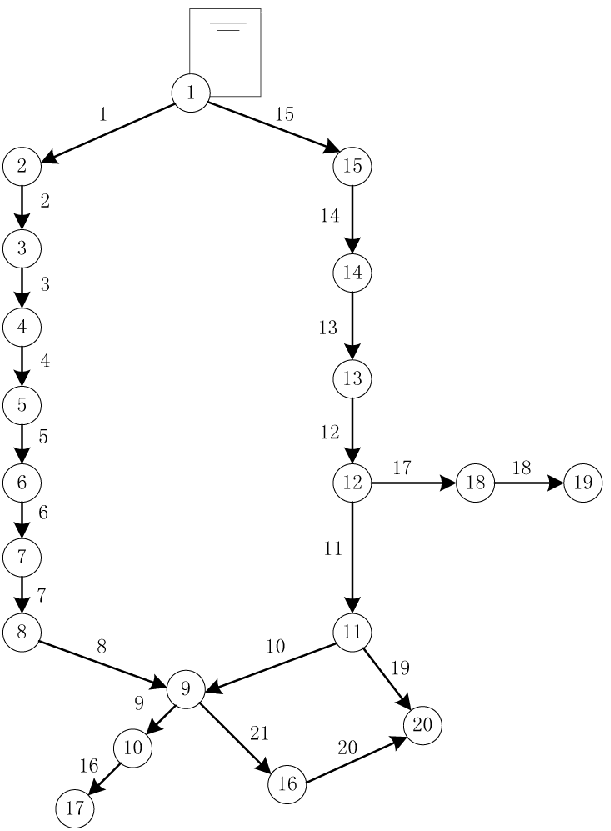}\\
  \caption{the New York network}\label{newyork}
\end{figure}

The layout of the New York network is given in Fig. \ref{newyork}. There are $20$ nodes, $21$ pipes and 1 loop in this network system.
At node 1, there exists a reservoir with 300-ft fixed head. The New York problem is different from other two cases, because there already exist pipes in the old system.
The common objective of this problem is to determine additional parallel pipes added to the existing ones to meet increased water demands while maintaining the minimum pressure requirements.
The the cost data, pipe and node data are given in Table \ref{costdatanewyork} and Table \ref{pipenodedatanewyork}, respectively. The Hazen-Williams coefficient $C$ is assumed to be $100$ in this case.

\begin{table}[!htbp]
\renewcommand{\arraystretch}{1.3}
\caption{Cost data for the New York network}
\label{costdatanewyork}
\centering
\footnotesize
\begin{tabular}{{p{0.8cm}p{1.2cm}p{1.2cm}p{0.8cm}p{1.2cm}p{1.2cm}}}
\hline
No. & Diameter (in.)   & Cost (\$/feet) & No. & Diameter (in.)   & Cost (\$/feet)\\
\hline
1 & 0  & 0.00 & 9 & 120 & 417.0\\
2 & 36 & 93.5 & 10& 132 & 469.0\\
3 & 48 & 134.0& 11& 144 & 522.0\\
4 & 60 & 176.0& 12& 156 & 577.0\\
5 & 72 & 221.0& 13& 168 & 632.0\\
6 & 84 & 267.0& 14& 180 & 689.0\\
7 & 96 & 316.0& 15& 192 & 746.0\\
8 & 108&365.0 & 16& 204 & 804.0\\
\hline
\end{tabular}
\end{table}

\begin{table*}[!htbp]
\centering
\begin{threeparttable}[b]
\renewcommand{\arraystretch}{1.3}
\caption{Pipe and node data for the New York network}
\label{pipenodedatanewyork}
\footnotesize
\begin{tabular}{{cccccc}}
\hline
Pipe  & Length (feet)\tnote{4}  & Existing Diameters (in.) & Node  & Demand (feet$^3$/s) \tnote{5}& Minimum Total Head (feet) \\
\hline
1 & 11600 & 180 & 1 & -2017.5   & 300.0\\
2 & 19800 & 180 & 2 & 92.4      & 255.0\\
3 & 7300  & 180 & 3 & 92.4      & 255.0\\
4 & 8300  & 180 & 4 & 88.2      & 255.0\\
5 & 8600  & 180 & 5 & 88.2      & 255.0\\
6 & 19100 & 180 & 6 & 88.2      & 255.0\\
7 & 9600  & 132 & 7 & 88.2      & 255.0\\
8 & 12500 & 132 & 8 & 88.2      & 255.0\\
9 & 9600  & 180 & 9 & 170.0     & 255.0\\
10& 11200 & 204 & 10& 1.0       & 255.0\\
11& 14500 & 204 & 11& 170.0     & 255.0\\
12& 12200 & 204 & 12& 117.1     & 255.0\\
13& 24100 & 204 & 13& 117.1     & 255.0\\
14& 21100 & 204 & 14& 92.4      & 255.0\\
15& 15500 & 204 & 15& 92.4      & 255.0\\
16& 26400 & 72  & 16& 170.0     & 260.0\\
17& 31200 & 72  & 17& 57.5      & 272.8\\
18& 24000 & 60  & 18& 117.1     & 255.0\\
19& 14400 & 60  & 19& 117.1     & 255.0\\
20& 38400 & 60  & 20& 170.0     & 255.0\\
21& 26400 & 72  & - &    -      & -    \\
\hline
\end{tabular}
\begin{tablenotes}
\item [4] 1 feet = 0.3048 m
\item [5] 1 feet$^3$/s = 28.3168 L/s
\end{tablenotes}
\end{threeparttable}
\end{table*}

\begin{figure*}[!htbp]
  \centering
  % Requires \usepackage{graphicx}
  \includegraphics[width=8cm]{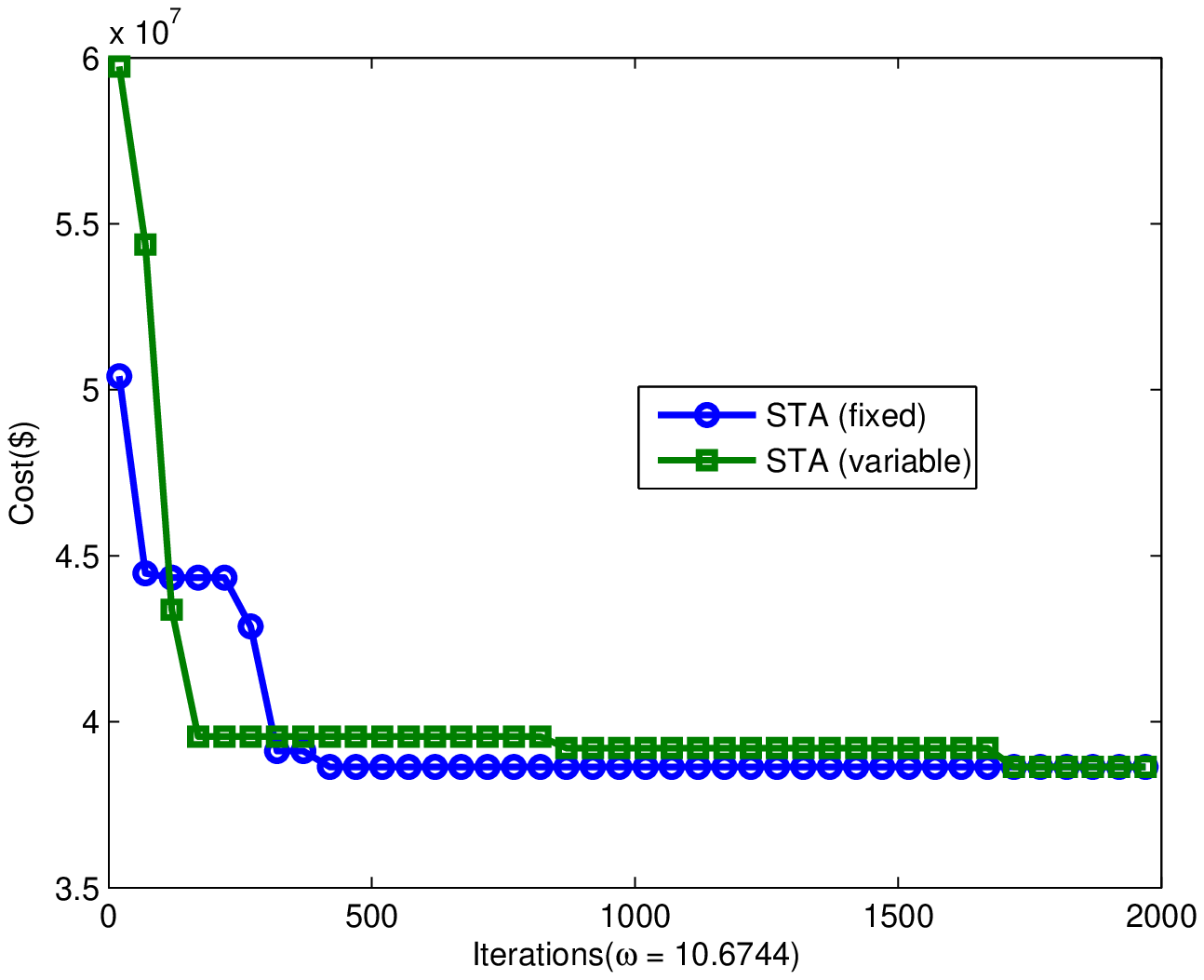}
  \includegraphics[width=8cm]{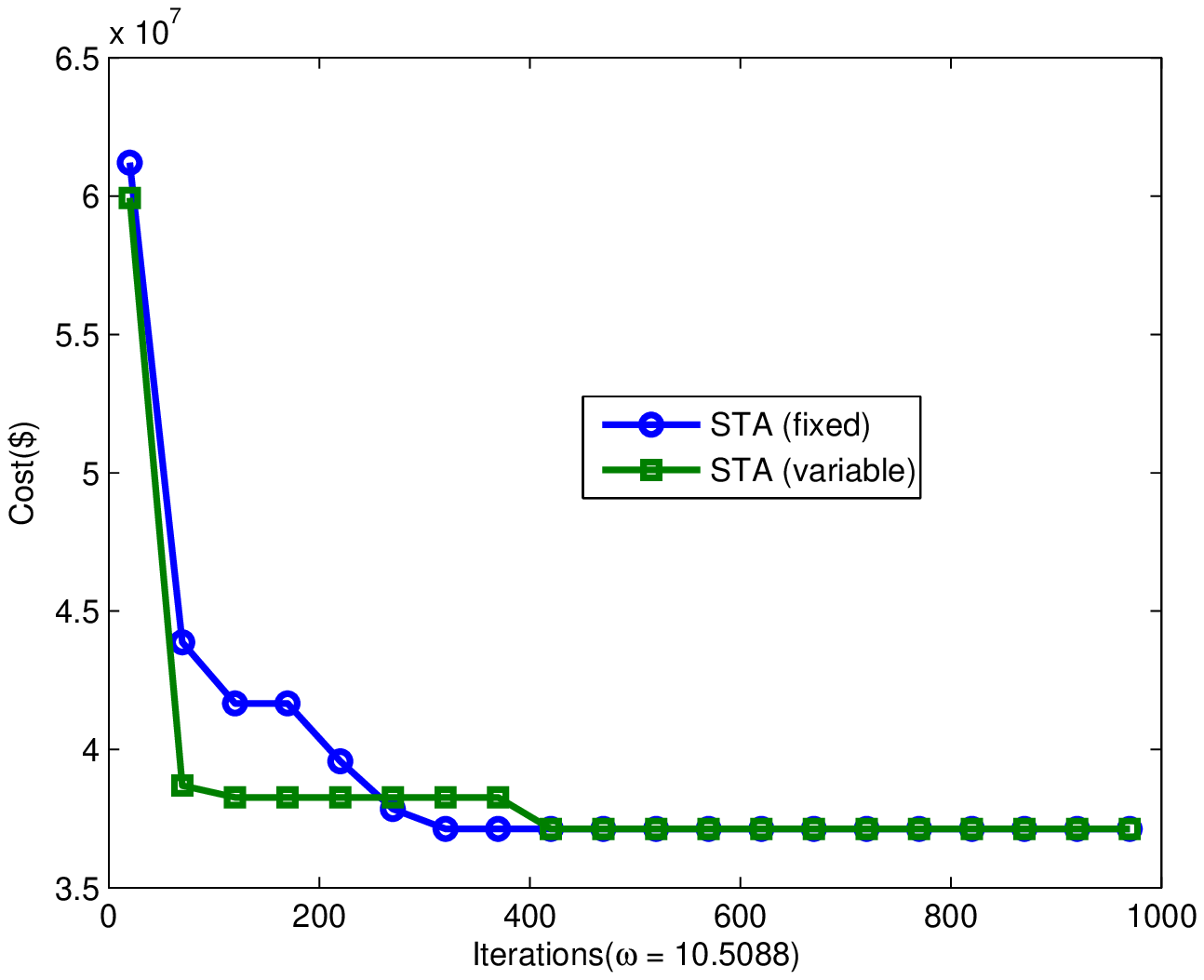}\\
  \caption{Iterative curves of two best solutions using STA for the New York problem when $\omega$ is 10.6744 and 10.5088, respectively}\label{newyorkcurveruns}
\end{figure*}

For the New York network, the search enforcement $SE$ is also set at $10$, and the penalty coefficient $pc$ is fixed at $2e6$, or varies from $1e6$ to $1e7$ in a linearly increasing way.
The maximum number of iterations is set at $2e3$, and a total of $20$ runs are executed for both fixed and variable $pc$. Fig. \ref{newyorkcurveruns} gives the
iterative curves of two best solutions and changes in $20$ runs for the New York problem with fixed and variable $pc$ respectively. We can find that the best solution is hit five times by the STA with fixed $pc$ and twice by the STA with variable $pc$.

\begin{remark}
Under the circumstance, the minimum function evaluations to achieve the best known solution is $5200$, which takes up $2.6883e$-$19\%$ of all possible combinations ($16^{21} = 1.9343e25$).
\end{remark}

Table \ref{newyorksolution} gives the best solutions gained by various algorithms, and it can be found that STA with both fixed and variable $pc$ can achieve the best known solution at the cost of $37.13$ million dollars. As a matter of fact, the same solution was also gained by GA \cite{savic1997} with the function evaluations at $1,000,000$. The pressure heads for the New York network obtained by the discrete STA are given in Table \ref{newyorkpressureheads}.

\begin{table*}[!htbp]
\caption{Solutions for the New York network}
\label{newyorksolution}
\renewcommand{\arraystretch}{1.3}
\centering
\footnotesize
\begin{tabular}{{p{2.0cm}p{1.2cm}p{1.5cm}p{1.2cm}p{1.8cm}p{1.8cm}p{1.8cm}p{1.8cm}}}
\hline
Pipe  & Gessler \cite{gessler1985} & Morgan and Goulter \cite{morgan1985} & Dandy ~\;\;et al. \cite{dandy1996} & STA (fixed) $\omega = 10.6744$& \;$\omega = 10.5088$ & STA (variable) $\omega = 10.6744$& \;$\omega = 10.5088$  \\
\hline
1  &   0 &   0 & 0   & 0   &  0  & 0   & 0 \\
2  &   0 &   0 & 0   & 0   &  0  & 0   & 0 \\
3  &   0 &   0 & 0   & 0   &  0  & 0   & 0 \\
4  &   0 &   0 & 0   & 0   &  0  & 0   & 0 \\
5  &   0 &   0 & 0   & 0   &  0  & 0   & 0 \\
6  &   0 &   0 & 0   & 0   &  0  & 0   & 0 \\
7  & 100 & 144 & 0   & 144 & 108 & 144 & 108 \\
8  & 100 & 144 & 0   & 0   & 0   & 0   & 0 \\
9  &   0 & 0   & 0   & 0   & 0   & 0   & 0 \\
10 &   0 & 0   & 0   & 0   & 0   & 0   & 0 \\
11 &   0 & 0   & 0   & 0   & 0   & 0   & 0 \\
12 &   0 & 0   & 0   & 0   & 0   & 0   & 0 \\
13 &   0 & 0   & 0   & 0   & 0   & 0   & 0 \\
14 &   0 & 0   & 0   & 0   & 0   & 0   & 0 \\
15 &   0 & 0   & 120 & 0   & 0   & 0   & 0 \\
16 & 100 & 96  & 84  & 96  & 96  & 96  & 96 \\
17 & 100 & 96  & 96  & 96  & 96  & 96  & 96 \\
18 & 80  & 84  & 84  & 84  & 96  & 84  & 84 \\
19 & 60  & 60  & 72  & 72  & 72  & 72  & 72 \\
20 &   0 & 0   & 0   & 0   & 0   & 0   & 0 \\
21 & 80  & 84  & 72  & 72  & 72  & 72  & 72 \\
Cost(\$ millions) & 41.80 & 39.20 & 38.80 & 38.64 &  37.13 & 38.64 & 37.13\\
\hline
\end{tabular}
\end{table*}

\begin{table*}[!htbp]
\caption{Pressure Heads for the New York network}
\label{newyorkpressureheads}
\renewcommand{\arraystretch}{1.3}
\centering
\footnotesize
\begin{tabular}{{p{1.0cm}p{1.8cm}p{1.8cm}p{1.8cm}p{1.8cm}p{1.8cm}p{1.8cm}p{1.8cm}}}
\hline
Node  & STA $\omega = 10.6744$& \;$\omega = 10.5088$ & Node &  STA $\omega = 10.6744$& \;$\omega = 10.5088$ \\
\hline
1  &  300.00 & 300.00  & 11  & 273.85 & 273.86\\
2  &  294.20 & 294.33  & 12  & 275.12 & 275.15\\
3  &  286.14 & 286.47  & 13  & 278.09 & 278.12\\
4  &  283.78 & 284.16  & 14  & 285.55 & 285.58\\
5  &  281.68 & 282.13  & 15  & 293.32 & 293.34\\
6  &  280.06 & 280.55  & 16  & 260.05 & 260.16\\
7  &  277.50 & 278.08  & 17  & 272.85 & 272.86\\
8  &  276.65 & 276.51  & 18  & 261.15 & 261.30\\
9  &  273.76 & 273.76  & 19  & 255.02 & 255.21\\
10 &  273.73 & 273.73  & 20  & 260.70 & 260.81\\
\hline
\end{tabular}
\end{table*}

\section{Conclusion}
The complexity of the water distribution network comes from two aspects, one is the linear and nonlinear equations, which are commonly handled by a hydraulic solver to ensure that the continuity and head loss equations are satisfied automatically, the other difficulty is that the commercial pipe size is discrete, which is proved to be NP-hard.

In this paper, It is shown that the network system can be reduced to the dimensionality of the number of closed simple loop or
required independent paths, which can reduce the computational complexity of solving linear and nonlinear equations simultaneously to a large extent.

To overcome the NP-hardness, a new intelligent optimization algorithm named discrete state transition algorithm is introduced to
find the optimal or suboptimal solution. There are four intelligent operators in discrete STA, which are easy to understand and to be implemented. The ``restore in probability"  $p_1$ and ``risk in probability" $p_2$ strategy in discrete STA is used to escape local optimal and increase the probability to capture the global optimum.

At first, a Monte Carlo simulation is studied to investigate a good combination of $p_1$ and $p_2$, and we find that $(p_1, p_2) = (0.1, 0.1)$ is a good choice. We then focus on a empirical study of the Two-Loop network, by training the network, we find that the penalty coefficient plays a significant role in the search ability and solution feasibility.

Based on the experience gained from the Two-Loop problem, the discrete STA has successfully applied to the Hanoi and New York networks, and the results show that the discrete STA can achieve the best known solutions with less function evaluations. The success of the discrete STA in optimal design of water distribution network has demonstrate that the discrete STA is a promising alternative in combinational optimization.
\section*{Acknowledgment}
The authors would like to thank...

% Can use something like this to put references on a page
% by themselves when using endfloat and the captionsoff option.
\ifCLASSOPTIONcaptionsoff
  \newpage
\fi


\begin{thebibliography}{1}

\bibitem{yates1984}
D. F. Yates, A. B. Templeman and T. B. Boffey, ``The computational complexity of the problem of determining least capital cost designs for water supply networks,"
\textit{Engineering optimization}, vol. 7, no. 2, pp. 142--155, 1984.


\bibitem{gessler1985}
J. Gessler, ``Pipe network optimization by enumeration,"
\textit{Computer applications in water resources}, ASCE, New York, pp. 572--581, 1985.


\bibitem{alperovits1977}
E. Alperovits and U. Shamir, ``Design of optimal water distribution systems,"  \textit{Water resources research}, vol. 13, no. 6, pp. 885--900, 1977.


\bibitem{morgan1985}
D. R. Morgan and I.C. Goulter,  ``Optimal urban water distribution design,"
\textit{Water Resources Research}, vol. 21, no. 5, pp. 642--652, 1985.


\bibitem{goulter1986}
I. C. Goulter, B. M. Lussier, and D. R. Morgan, ``Implications of head loss path choice in the optimization", \textit{Water Resources Research}, vol. 22, no. 5, pp. 819--822, 1986.


\bibitem{kessler1989}
A. Kessler and U. Shamir, ``Analysis of the linear programming gradient method for optimal design of water supply networks", \textit{Water Resources Research}, vol. 25, no. 7, pp. 1469--1480, 1989.


\bibitem{eiger1994}
G. Eiger, U. Shamir, and A. Ben-Tal,  ``Optimal design of water distribution networks", \textit{Water resources research}, vol. 30, no. 9, pp. 2637--2646, 1994.


\bibitem{dandy1996}
G. C. Dandy, A. R. Simpson and L. J. Murphy,  ``An improved genetic algorithm for pipe network optimization", \textit{Water Resources Research}, vol. 32, no. 2, pp. 449--458, 1996.

\bibitem{savic1997}
D. A. Savic and G. A. Walters, ``Genetic algorithms for least-cost design of water distribution networks", \textit{Journal of Water Resources Planning and Management}, vol. 123, no. 2, pp. 67--77, 1997.


\bibitem{simpson1994}
A. R. Simpson, G. C.  Dandy and L. J. Murphy, ``Genetic algorithms compared to other techniques for pipe optimization", \textit{Journal of Water Resources Planning and Management}, vol. 120, no. 4, pp. 423--443, 1994.

\bibitem{cunha1999}
M. da Concei{\c{c}}{\~a}o Cunha and J. Sousa, ``Water distribution network design optimization: Simulated annealing approach", \textit{Journal of Water Resources Planning and Management}, vol. 125, no. 4, pp. 215--221, 1999.


\bibitem{liong2004}
S. Y. Liong and M. Atiquzzaman, ``Optimal design of water distribution network using shuffled complex evolution", \textit{Journal of The Institution of Engineers, Singapore}, vol. 44, no. 1, pp. 93--107, 2004.

\bibitem{zecchin2005}
A. C. Zecchin, A. R. Simpson, H. R. Maier and J. B. Nixon, ``Parametric study for an ant algorithm applied to water distribution system optimization", \textit{IEEE Transactions on Evolutionary Computation}, vol. 9, no. 2, pp. 175--191, 2005.


\bibitem{zecchin2006}
A. C. Zecchin, A.R. Simpson, H. R. Maier, M. Leonard, A. J. Roberts and M. J. Berrisford,  ``Application of two ant colony optimisation algorithms to water distribution system optimisation", \textit{Mathematical and computer modelling}, vol. 44, no. 5, pp. 451--468, 2006.


\bibitem{geem2006}
Z. W. Geem, ``Optimal cost design of water distribution networks using harmony search", \textit{Engineering Optimization}, vol. 38, no. 3, pp. 259--277, 2006.


\bibitem{montalvo2008}
I. Montalvo, J. Izquierdo, R. P{\'e}rez, and M. M. Tung, ``Particle Swarm Optimization applied to the design of water supply systems", \textit{Computers \& Mathematics with Applications}, vol. 56, no. 3, pp. 769--776, 2008.

\bibitem{vasan2010}
A. Vasan and S. P. Simonovic, ``Optimization of water distribution network design using differential evolution", \textit{Journal of Water Resources Planning and Management}, vol. 136, no. 2, pp. 279--287, 2010.

\bibitem{cisty2010}
M. Cisty, ``Hybrid genetic algorithm and linear programming method for least-cost design of water distribution systems", \textit{Water resources management}, vol. 24, no. 1, pp. 1--24, 2010.

\bibitem{haghighi2011}
A. Haghighi, H. M. V. Samani and Z. M. V. Samani, ``GA-ILP method for optimization of water distribution networks", \textit{Water resources management}, vol. 25, no. 7, pp. 1791--1808, 2011.

\bibitem{xzhou2011a}
X. J. Zhou,  C. H. Yang, W. H. Gui, ``Initial version of state transition algorithm", in International conference on digital manufacturing and automation (ICDMA), 644--647, 2011.

\bibitem{xzhou2011b}
X. J. Zhou,  C. H. Yang, W. H. Gui, ``A new transformation into state transition algorithm for finding the global minimum", in International conference on intelligent control and
information processing (ICICIP), 674--678, 2011.

\bibitem{xzhou2012}
X. J. Zhou,  C. H. Yang, W. H. Gui, ``State transition algorithm", Journal of Industrial and Management Optimization, vol. 8, no. 4, pp. 1039--1056, 2012.

\bibitem{yang2012}
C. H. Yang, X. L. Tang, X. J. Zhou and W. H. Gui, ``State transition algorithm for traveling salesman problem", in Proceeding of the 31st Chinese Control Conference, pp. 2481--2485, 2012.

\bibitem{metroplis1953}
Metropolis, N., Rosenbluth, A., Rosenbluth, M., Teller, A., Teller, E., Equation of state calculations by fast computing machines, Journal of Chemical Physics, 21(6), 1087--1092, 1953.

\bibitem{andradottir1996}
Andrad{\'o}ttir, S., A global search method for discrete stochastic optimization, SIAM Journal on Optimization, 6(2), 513--530, 1996.

\bibitem{andradottir1999}
Andrad{\'o}ttir, S., Accelerating the convergence of random search methods for discrete stochastic optimization, ACM Transactions on Modeling and Computer Simulation, 9(4), 349--380, 1999.

\bibitem{baba}
Baba, N., Shoman T., Sawaragi Y., A modified convergence theorem for a random optimization method, Information Sciences, 13(2), 159--166, 1977.

\bibitem{van}
Van Laarhoven P.J.M.,  Aarts E.H.L., Simulated annealing: theory and applications, Dordrecht, The Netherlands, D. Reidel, 1987.

\bibitem{rossman1994}
Rossman, L.A., EPANET 2: users manual, US Environmental Protection Agency, Office of Research and Development, National Risk Management Research Laboratory, 1994.

\bibitem{wood1980}
Wood, D.J., Computer Analysis of Flow in Pipe Networks Including Extended Period Simulations: User's Manual, Office of Continuing Education and Extension of the College of Engineering of the University of Kentucky, 1980.




\end{thebibliography}
\end{document}